\pdfoutput=1
\documentclass[11pt,a4paper]{article}
\usepackage[latin1]{inputenc}
\usepackage{eucal,amssymb,amsfonts,amsmath,amsthm,graphicx,bbm,url}
\usepackage{hyperref}
\usepackage{fullpage}

\newtheorem{Theorem}{Theorem}

\newtheorem{Remark}[Theorem]{Remark}
\newtheorem{Lemma}[Theorem]{Lemma}

\newtheorem{Proposition}[Theorem]{Proposition}

\newtheorem{Open question}[Theorem]{Open question}

\begin{document}

\title{On the critical value function in the divide and color model}

\author{Andr{\'a}s B{\'a}lint\thanks{Chalmers University of Technology,
    e-mail: \url{abalint@few.vu.nl}} \and
  Vincent Beffara\thanks{UMPA-ENS Lyon,
    e-mail: \url{vbeffara@ens-lyon.fr}} \and
  Vincent Tassion\thanks{ENS Lyon,
    e-mail: \url{vincent.tassion@ens-lyon.fr}}}

\maketitle

\begin{abstract}
  The divide and color model on a graph $G$ arises by first deleting
  each edge of $G$ with probability $1-p$ independently of each other,
  then coloring the resulting connected components (\emph{i.e.}, every
  vertex in the component) black or white with respective probabilities
  $r$ and $1-r$, independently for different components.  Viewing it as
  a (dependent) site percolation model, one can define the critical
  point $r_c^G(p)$.

  In this paper, we mainly study the continuity properties of the
  function $r_c^G$, which is an instance of the question of locality
  for percolation. Our main result is the fact that in the case
  $G=\mathbb Z^2$, $r_c^G$ is continuous on the interval $[0,1/2)$; we
  also prove continuity at $p=0$ for the more general class of graphs
  with bounded degree. We then investigate the sharpness of the
  bounded degree condition and the monotonicity of $r_c^G(p)$ as a
  function of $p$.
\end{abstract}




\section*{Introduction}

The divide and color (DaC) model is a natural dependent site
percolation model introduced by H{\"a}ggstr{\"o}m in
\cite{HaggstromDaC}.  It has been studied directly in
\cite{HaggstromDaC,Garet,BCM,BBT}, and as a member of a more general
family of models in \cite{KW,BCM,BGibbs,GrGr}.  This model is defined
on a multigraph $G=(\mathcal{V},\mathcal{E})$, where $\mathcal{E}$ is
a multiset (\emph{i.e.}, it may contain an element more than once),
thus allowing parallel edges between pairs of vertices.  For
simplicity, we will imprecisely call $G$ a \emph{graph} and
$\mathcal{E}$ the \emph{edge set}, even if $G$ contains self-loops or
multiple edges.  The DaC model with parameters $p,r \in [0,1]$, on a
general (finite or infinite) graph $G$ with vertex set $\mathcal{V}$
and edge set $\mathcal{E}$, is defined by the following two-step
procedure:

\begin{itemize}
\item First step: Bernoulli bond percolation.  We independently declare
  each edge in $\mathcal{E}$ to be open with probability $p$, and closed
  with probability $1-p$. We can identify a bond percolation
  configuration with an element $\eta \in \{0,1\}^{\mathcal{E}}$: for
  each $e\in \mathcal{E}$, we define $\eta (e)=1$ if $e$ is open, and
  $\eta (e)=0$ if $e$ is closed.
\item Second step: Bernoulli site percolation on the resulting cluster
  set.  Given $\eta \in \{0,1\}^{\mathcal{E}}$, we call
  \emph{$p$-clusters} or \emph{bond clusters} the connected components
  in the graph with vertex set $\mathcal{V}$ and edge set $\{e\in
  \mathcal{E}:\eta (e)=1\}$.  The set of $p$-clusters of $\eta$ gives
  a partition of $\mathcal{V}$.  For each $p$-cluster $\mathcal{C}$,
  we assign the same color to all the vertices in $\mathcal{C}$. The
  chosen color is black with probability $r$ and white with
  probability $1-r$, and this choice is independent for different
  $p$-clusters.
\end{itemize}

These two steps yield a site percolation configuration $\xi \in
\{0,1\}^{\mathcal{V}}$ by defining, for each $v\in \mathcal{V}$, $\xi
(v)=1$ if $v$ is black, and $\xi (v)=0$ if $v$ is white.  The connected
components (via the edge set $\mathcal{E}$) in $\xi$ of the same color
are called (black or white) \emph{$r$-clusters}. The resulting measure
on $\{0,1\}^{\mathcal{V}}$ is denoted by $\mu_{p,r}^{G}$.

\bigskip

Let $E_{\infty }^b\subset \{0,1\}^{\mathcal{V}}$ denote the event that
there exists an infinite black $r$-cluster. By standard arguments (see
Proposition 2.5 in \cite{HaggstromDaC}), for each $p\in [0,1]$, there
exists a \emph{critical coloring value} $r_c^G(p)\in [0,1]$ such that
\[\mu _{p,r}^G(E_{\infty }^b) \begin{cases}
  =0 &\textrm{if } r<r_c^G(p),\\
  >0 &\textrm{if } r>r_c^G(p).\\
\end{cases}\]

The \emph{critical edge parameter} $p_c^G\in [0,1]$ is defined as
follows: the probability that there exists an infinite bond cluster is
$0$ for all $p<p_c^G$, and positive for all $p>p_c^G$. The latter
probability is in fact $1$ for all $p>p_c^G$, whence $r_c^G(p)=0$ for
all such $p$. Kolmogorov's $0-1$ law shows that in the case when all
the bond clusters are finite, $\mu _{p,r}^G(E_{\infty }^b)\in
\{0,1\}$; nevertheless it is possible that $\mu_{p,r}^G(E_{\infty
}^b)\in (0,1)$ for some $r>r_c^G(p)$ (\emph{e.g.} on the square
lattice, as soon as $p>p_c=1/2$, one has $\mu_{p,r}^G(E_{\infty }^b) =
r$).

\subsection*{Statement of the results}

Our main goal in this paper is to understand how the critical coloring
parameter $r_c^G$ depends on the edge parameter $p$.  Since the
addition or removal of self-loops obviously does not affect the value
of $r_c^G(p)$, we will assume that all the graphs $G$ that we consider
are without self-loops.  On the other hand, $G$ is allowed to contain
multiple edges.

Our first result, based on a stochastic domination argument, gives
bounds on $r_c^G(p)$ in terms of $r_c^G(0)$, which is simply the
critical value for Bernoulli site percolation on $G$. By the
\emph{degree} of a vertex $v$, we mean the number of edges incident on
$v$ (counted with multiplicity).

\begin{Proposition}\label{rcbounds}
  For any graph $G$ with maximal degree $\Delta $, for all $p\in
  [0,1)$,
  \[ 1 - \frac{1-r_c^G(0)}{(1-p)^{\Delta}} \leq r_c^G(p) \leq
  \frac{r_c^G(0)}{(1-p)^{\Delta}}. \]
\end{Proposition}

As a direct consequence, we get continuity at $p=0$ of the critical
value function:

\begin{Proposition}\label{contin0}
  For any graph $G$ with bounded degree, $r_c^G(p)$ is continuous in
  $p$ at $0$.
\end{Proposition}

One could think of an alternative approach to the question, as
follows: the DaC model can be seen as Bernoulli site percolation of
the random graph $G_p=(V_p,E_p)$ where $V_p$ is the set of
bond clusters and two bond clusters are connected by a bond of $E_p$
if and only if they are adjacent in the original graph. The study of
how $r_c^G(p)$ depends on $p$ is then a particular case of a more
general question known as the \emph{locality problem}: is it true in
general that the critical points of site percolation on a graph and a
small perturbation of it are always close? Here, for small $p$, the
graphs $G$ and $G_p$ are somehow very similar, and their critical
points are indeed close.

Dropping the bounded-degree assumption allows for the easy
construction of graphs for which continuity does not hold at $p=0$:

\begin{Proposition}\label{nonbounded}
  There exists a graph $G$ with $p_c^G>0$ such that $r_c^G$ is
  discontinuous at $0$.
\end{Proposition}

In general, when $p>0$, the graph $G_p$ does not have bounded degree,
even if $G$ does; this simple remark can be exploited to construct
bounded degree graphs for which $r_c^G$ has discontinuities below the
critical point of bond percolation (though of course not at $0$):

\begin{Theorem}\label{boundeddegreeconstruction}
  There exists a graph $G$ of bounded degree satisfying $p_c^G > 1/2$
  and such that $r_c^G(p)$ is discontinuous at $1/2$.
\end{Theorem}

\begin{Remark}
  The value $1/2$ in the statement above is not special: in fact, for
  every $p_0 \in (0,1)$, it is possible to generalize our argument to
  construct a graph with a critical bond parameter above $p_0$ and for
  which the discontinuity of $r_c$ occurs at $p_0$.
\end{Remark}

Our main results concerns the case $G=\mathbb Z^2$, for which the
above does not occur:

\begin{Theorem}\label{continp}
  The critical coloring value $r_c^{\mathbb{Z}^2}(p)$ is a continuous
  function of $p$ on the whole interval $[0,1/2)$.
\end{Theorem}

The other, perhaps more anecdotal question we investigate here is
whether $r_c^G$ is monotonic below $p_c$. This is the case on the
triangular lattice (because it is constant equal to $1/2$), and
appears to hold on $\mathbb Z^2$ in simulations (see the companion
paper \cite{BBT}).

In the general case, the question seems to be rather
delicate. Intuitively the presence of open edges would seem to make
percolation easier, leading to the intuition that the function $p
\mapsto r_c(p)$ should be nonincreasing. Theorem 2.9 in
\cite{HaggstromDaC} gives a counterexample to this intuition. It is
even possible to construct quasi-transitive graphs on which any
monotonicity fails:

\begin{Proposition}\label{PropNonMononicity}
  There exists a quasi-transitive graph $G$ such that $r_c^G$ is not
  monotone on the interval $[0,p_c^{G})$.
\end{Proposition}

A brief outline of the paper is as follows.  We set the notation and
collect a few results from the literature in
Section~\ref{definitions-section}.  In
Section~\ref{section-stochastic-domination}, we stochastically compare
$\mu _{p,r}^{G}$ with Bernoulli site percolation
(Theorem~\ref{stdom}), and show how this result implies
Proposition~\ref{rcbounds}. We then turn to the proof of
Theorem~\ref{continp} in Section~\ref{contz2-section}, based on a
finite-size argument and the continuity of the probability of
cylindrical events.

In Section~\ref{tree-like-section}, we determine the critical value
function for a class of tree-like graphs, and in the following section
we apply this to construct most of the examples of graphs we mentioned
above.

\section{Definitions and notation}\label{definitions-section}

We start by explicitly constructing the model, in a way which will be
more technically convenient than the intuitive one given in the
introduction.

Let $G$ be a connected graph $( \mathcal{V},\mathcal{E} )$ where the
set of vertices $\mathcal{V}=\{v_0,v_1,v_2,\ldots\}$ is countable.  We
define a total order ``$<$'' on $\mathcal{V}$ by saying that $v_i<v_j$
if and only if $i<j$.  In this way, for any subset $V \subset
\mathcal{V}$, we can uniquely define $\min (V)\in V$ as the minimal
vertex in $V$ with respect to the relation ``$<$''.  For a set $S$, we
denote $\{0,1\}^S$ by $\Omega _S$. We call the elements of $\Omega
_{\mathcal {E}}$ {\em bond configurations}, and the elements of
$\Omega _{\mathcal {V}}$ \emph{site configurations}.  As defined in
the Introduction, in a bond configuration $\eta $, an edge $e\in
\mathcal{E}$ is called \emph{open} if $\eta (e)=1$, and \emph{closed}
otherwise; in a site configuration $\xi $, a vertex $v\in \mathcal{V}$
is called \emph{black} if $\xi (e)=1$, and \emph{white} otherwise.
Finally, for $\eta \in \Omega _{\mathcal {E}}$ and $v\in \mathcal{V}$,
we define the \emph{bond cluster} $\mathcal{C}_v(\eta )$ of $v$ as the
maximal connected induced subgraph containing $v$ of the graph with
vertex set $\mathcal{V}$ and edge set $\{e\in \mathcal{E}:\eta
(e)=1\}$, and denote the vertex set of $\mathcal{C}_v(\eta )$ by
$C_v(\eta )$.

For $a\in [0,1]$ and a set $S$, we define $\nu _a ^S$ as the probability
measure on $\Omega _S$ that assigns to each $s\in S$ value $1$ with
probability $a$ and $0$ with probability $1-a$, independently for
different elements of $S$.  We define a function
\begin{equation*}
  \begin{array}{lccc}
    \Phi \: : \: &\Omega_\mathcal{E} \times \Omega_\mathcal{V}
    &\rightarrow& \Omega_\mathcal{E} \times \Omega_\mathcal{V},\\
    &(\eta,\kappa)&\mapsto&(\eta,\xi),\\
\end{array}
\end{equation*}
where $\xi(v)=\kappa(\min(C_v(\eta )))$.  For $p,r\in [0,1]$, we define
$\mathbb{P}^G_{p,r}$ to be the image measure of $\nu_p^\mathcal{E}
\otimes \nu_r^\mathcal{V}$ by the function $\Phi $, and denote by $\mu
^G_{p,r}$ the marginal of $\mathbb{P}^G_{p,r}$ on $\Omega
_{\mathcal{V}}$.  Note that this definition of $\mu ^G_{p,r}$ is
consistent with the one in the Introduction.

Finally, we give a few definitions and results that are necessary for
the analysis of the DaC model on the square lattice, that is the graph
with vertex set $\mathbb{Z}^2$ and edge set $\mathcal{E}^2=\{\left<
  v,w\right> :v=(v_1,v_2),w=(w_1,w_2)\in \mathbb{Z}^2,\
|v_1-w_1|+|v_2-w_2|=1\}$. The \emph{matching graph}
$\mathbb{Z}^2_*$ {of the square lattice} is the graph with vertex set
$\mathbb{Z}^2$ and edge set $\mathcal{E}^2_*=\{\left< v,w\right>
:v=(v_1,v_2),w=(w_1,w_2)\in \mathbb{Z}^2,\ \max
(|v_1-w_1|,|v_2-w_2|)=1\}$. In the same manner as in the Introduction,
we define, for a color configuration $\xi \in \{0,1\}^{\mathbb{Z}^2}$,
(black or white) \emph{$*$-clusters} as connected components (via the
edge set $\mathcal{E}^2_*$) in $\xi $ of the same color. We denote by
$\Theta ^*(p,r)$ the $\mathbb{P}_{p,r}^{\mathbb{Z}^2}$-probability
that the origin is contained in an infinite black $*$-cluster, and
define \[r_c^*(p) = \sup \{r:\Theta ^*(p,r)=0\}\] for all $p \in
[0,1]$ --- note that this value may differ from
$r_c^{\mathbb{Z}^2_*}(p)$. The main result in \cite{BCM} is that for
all $p\in [0,1/2)$, the critical values $r_c^{\mathbb{Z}^2}(p)$ and
$r_c^*(p)$ satisfy the duality relation
\begin{equation}\label{rcplusrcstaris1}
  r_c^{\mathbb{Z}^2}(p) + r_c^*(p)= 1.
\end{equation}

We will also use exponential decay result for subcritical Bernoulli
bond percolation on $\mathbb{Z}^2$.  Let ${\bf 0}$ denote the origin
in $\mathbb{Z}^2$, and for each $n\in \mathbb{N}=\{1,2,\ldots \}$, let
us define $S_n=\{v\in \mathbb{Z}^2:dist(v,{\bf 0})=n\}$ (where $dist$
denotes graph distance), and the event $M_n=\{\eta \in \Omega
_{\mathcal{E}^2}:$ there is a path of open edges in $\eta $ from ${\bf
  0}$ to $S_n\}$. Then we have the following result:

\begin{Theorem}[\cite{kesten}]\label{lexpdecay}
  For $p<1/2$, there exists $\psi (p)>0$ such that for all $n\in
  \mathbb{N}$, we have that
  \begin{displaymath}
    \nu _{p}^{\mathcal{E}^2}(M_n)< e^{-n\psi (p)}.
  \end{displaymath}
\end{Theorem}

\section{Stochastic domination and continuity at \texorpdfstring{$p=0$}{p=0}}
\label{section-stochastic-domination}

In this section, we prove Proposition~\ref{rcbounds} via a stochastic
comparison between the DaC measure and Bernoulli site percolation.
Before stating the corresponding result, however, let us recall the
concept of stochastic domination.

We define a natural partial order on $\Omega _{\mathcal{V}}$ by saying
that $\xi ^{\prime }\geq \xi $ for $\xi ,\xi ^{\prime }\in \Omega
_{\mathcal{V}}$ if, for all $v\in \mathcal{V}$, $\xi ^{\prime }(v)\geq
\xi (v)$.  A random variable $f:\Omega _{\mathcal{V}}\to \mathbb{R}$ is
called \emph{increasing} if $\xi ^{\prime }\geq \xi $ implies that
$f(\xi ^{\prime })\geq f(\xi )$, and an event $E\subset \Omega
_{\mathcal{V}}$ is increasing if its indicator random variable is
increasing.  For probability measures $\mu ,\mu ^{\prime }$ on $\Omega
_{\mathcal{V}}$, we say that $\mu ^{\prime }$ is \emph{stochastically
  larger} than $\mu $ (or, equivalently, that $\mu $ is
\emph{stochastically smaller} than $\mu ^{\prime }$, denoted by $\mu
\leq _{\textrm{st}}\mu ^{\prime }$) if, for all bounded increasing
random variables $f:\Omega _{\mathcal{V}}\to \mathbb{R}$, we have that
$$
\int _{\Omega _{\mathcal{V}}}f(\xi )\ d\mu ^{\prime }(\xi )\geq \int
_{\Omega _{\mathcal{V}}}f(\xi )\ d\mu (\xi ).
$$
By Strassen's theorem \cite{Strassen}, this is equivalent to the
existence of an appropriate coupling of the measures $\mu ^{\prime }$
and $\mu $; that is, the existence of a probability measure $\mathbb{Q}$
on $\Omega _{\mathcal{V}}\times \Omega _{\mathcal{V}}$ such that the
marginals of $\mathbb{Q}$ on the first and second coordinates are $\mu
^{\prime }$ and $\mu $ respectively, and $\mathbb{Q}(\{(\xi ^{\prime
},\xi )\in \Omega _{\mathcal{V}}\times \Omega _{\mathcal{V}}: \xi
^{\prime }\geq \xi \})=1$.

\begin{Theorem}\label{stdom}
  For any graph $G=( \mathcal{V},\mathcal{E} )$ whose maximal degree is
  $\Delta $, at arbitrary values of the parameters $p,r\in [0,1]$,
  \begin{displaymath}
    \nu _{r(1-p)^{\Delta }}^{\mathcal{V}} \leq _{\textrm{st}}\mu
    _{p,r}^G \leq _{\textrm{st}}\nu _{1-(1-r)(1-p)^{\Delta
      }}^{\mathcal{V}}.
  \end{displaymath}
\end{Theorem}

Before turning to the proof, we show how Theorem~\ref{stdom} implies
Proposition~\ref{rcbounds}.

\paragraph{Proof of Proposition~\ref{rcbounds}.}

It follows from Theorem~\ref{stdom} and the definition of stochastic
domination that for the increasing event $E_{\infty }^b$ (which was
defined in the Introduction), we have $\mu ^{G} _{p,r}(E_{\infty }^b)>0$
whenever $r(1-p)^{\Delta }>r_c^{G}(0)$, which implies that $r_c^G(p)
\leq r_c^G(0)/(1-p)^{\Delta }$.  The derivation of the lower bound for
$r_c^{G}(p)$ is analogous.  \qed

\bigskip

Now we give the proof of Theorem~\ref{stdom}, which bears some resemblance with
the proof of Theorem 2.3 in \cite{HaggstromDaC}.

\paragraph{Proof of Theorem~\ref{stdom}.}

Fix $G=(\mathcal{V},\mathcal{E})$ with maximal degree $\Delta $, and
parameter values $p,r\in [0,1]$.  We will use the relation ``$<$'' and
the minimum of a vertex set with respect to this relation as defined in
Section~\ref{definitions-section}.  In what follows, we will define
several random variables; we will denote the joint distribution of all
these variables by $\mathbb{P}$.

First, we define a collection $(\eta _{x,y}^e:x,y\in \mathcal{V},
e=\left< x,y\right> \in \mathcal{E})$ of i.i.d.\ Bernoulli($p$) random
variables (\emph{i.e.}, they take value $1$ with probability $p$, and
$0$ otherwise); one may imagine having each edge $e\in \mathcal{E}$
replaced by two directed edges, and the random variables represent
which of these edges are open. We define also a set $(\kappa _{x}:x\in
\mathcal{V})$ of Bernoulli($r$) random variables. Given a realization
of $(\eta _{x,y}^e:x,y\in \mathcal{V},e=\left< x,y\right> \in
\mathcal{E})$ and $(\kappa _{x}:x\in \mathcal{V})$, we will define an
$\Omega _{\mathcal{V}}\times \Omega _{\mathcal{E}}$-valued random
configuration $(\eta ,\xi )$ with distribution $\mathbb{P}
_{p,r}^{G}$, by the following algorithm.

\begin{enumerate}
\item Let $v=\min \{x\in \mathcal{V}:$ no $\xi $-value has been assigned
  yet to $x$ by this algorithm$\}$.  (Note that $v$ and $V, v_i, H_i $
  $(i\in \mathbb{N})$, defined below, are running variables,
  \emph{i.e.}, their values will be redefined in the course of the
  algorithm.)
\item We explore the ``directed open cluster'' $V$ of $v$ iteratively,
  as follows. Define $v_0=v$.  Given $v_0,v_1, \ldots ,v_i$ for some
  integer $i\geq 0$, set $\eta (e)=\eta ^e _{v_i,w}$ for every edge
  $e=\left< v_i,w\right> \in \mathcal{E}$ incident to $v_i$ such that no
  $\eta $-value has been assigned yet to $e$ by the algorithm, and write
  $H_{i+1}=\{w\in \mathcal{V}\setminus \{v_0,v_1,\ldots ,v_i\}:w$ can be
  reached from any of $v_0,v_1,\ldots ,v_i$ by using only those
  edges $e\in \mathcal{E}$ such that $\eta (e)=1$ has been assigned to
  $e$ by this algorithm$\}$.  If $H_{i+1}\neq \emptyset $, then we
  define $v_{i+1}=\min (H_{i+1})$, and continue exploring the directed
  open cluster of $v$; otherwise, we define $V=\{v_0,v_1,\ldots ,v_i\}$, and move to
  step 3.
\item Define $\xi (w)=\kappa _v$ for all $w\in V$, and return to step 1.
\end{enumerate}

It is immediately clear that the above algorithm eventually assigns a
$\xi$-value to each vertex. Note also that a vertex $v$ can receive a
$\xi$-value only after all edges incident to $v$ have already been
assigned an $\eta $-value, which shows that the algorithm eventually
determines the full edge configuration as well. It is easy to convince
oneself that $(\eta ,\xi )$ obtained this way indeed has the desired
distribution.

Now, for each $v\in \mathcal{V}$, we define $Z(v)=1$ if $\kappa _v=1$
and $\eta ^e _{w,v}=0$ for all edges $e=\left< v,w\right> \in \mathcal{E}$
incident on $v$ (\emph{i.e.}, all directed
edges towards $v$ are closed), and $Z(v)=0$ otherwise.  Note that every
vertex with $Z(v)=1$ has $\xi (v)=1$ as well, whence the distribution of
$\xi $ (\emph{i.e.}, $\mu _{p,r}^{G}$) stochastically dominates the
distribution of $Z$ (as witnessed by the coupling $\mathbb{P}$).

Notice that $Z(v)$ depends only on the states of the edges pointing to
$v$ and on the value of $\kappa_v$; in particular the distribution of
$Z$ is a product measure on $\Omega _{\mathcal{V}}$ with parameter
$r(1-p)^{d(v)}$ at $v$, where $d(v)\leq \Delta $ is the degree of $v$,
whence $\mu _{p,r}^{G}$ stochastically dominates the product measure
on $\Omega _{\mathcal{V}}$ with parameter $r(1-p)^{\Delta }$, which
gives the desired stochastic lower bound.  The upper bound can be
proved analogously; alternatively, it follows from the lower bound by
exchanging the roles of black and white. \qed

\section{Continuity of \texorpdfstring{$r_c^{\mathbb{Z}^2}(p)$}{rcp}
  on the interval
  \texorpdfstring{$[0,1/2)$}{[0,1/2)}}\label{contz2-section}

In this section, we will prove Theorem~\ref{continp}. Our first task
is to prove a technical result valid on more general graphs stating
that the probability of any event $A$ whose occurrence depends on a
finite set of $\xi$-variables is a continuous function of $p$ for
$p<p_c^G$. The proof relies on the fact that although the color of a
vertex $v$ may be influenced by edges arbitrarily far away, if
$p<p_c^G$, the corresponding influence decreases to $0$ in the limit
as we move away from $v$. Therefore, the occurrence of the event $A$
depends essentially on a finite number of $\eta $- and $\kappa
$-variables, whence its probability can be approximated up to an
arbitrarily small error by a polynomial in $p$ and $r$.

Once we have proved Proposition~\ref{cylcont} below, which is valid on
general graphs, we will apply it on $\mathbb{Z}^2$ to certain
``box-crossing events,'' and appeal to results in \cite{BCM} to deduce
the continuity of $r_c^{\mathbb{Z}^2}(p)$.

\begin{Proposition}\label{cylcont}
  For every site percolation event $A \subset \{0,1\}^{\mathcal{V}}$
  depending on the color of finitely many vertices, $\mu_{p,r}^G(A)$
  is a continuous function of $(p,r)$ on the set $[0,p_c^G) \times
  [0,1]$.
\end{Proposition}

\paragraph{Proof.}

In this proof, when $\mu$ is a measure on a set $S$, $X$ is a random
variable with law $\mu$ and $F:\: S\longrightarrow \mathbb{R}$ is a
bounded measurable function, we write
abusively $\mu[F(X)]$ for the expectation of $F(X)$.
We show a slightly more general result: for any
$k \geq 1$, $\boldsymbol{x}=(x_1,\ldots,x_k)\in \mathcal{V}^k$ and $f:
\{0,1\}^k \to \mathbb{R}$ bounded and measurable, $\mu_{p,r} ^G\left [
  f(\xi(x_1),\ldots,\xi(x_k))\right ]$ is continuous in $(p,r)$ on the
product $[0,p_c^G) \times [0,1]$.  Proposition~\ref{cylcont} will follow
by choosing an appropriate family $\{x_1,\ldots,x_k\}$ such that the
states of the $x_i$ suffices to determine whether $A$ occurs, and take
$f$ to be the indicator function of $A$.

To show the previous affirmation, we condition on the vector
\[
\boldsymbol{m}_{\boldsymbol{x}}(\eta)=(\min C_{x_1}(\eta),\ldots, \min
C_{x_k}(\eta))
\] which takes values in the finite set $\boldsymbol{V}=\left
  \{(v_1,\ldots,v_k) \in \mathcal{V}^k \: : \: \forall i \: v_i \leq
  \max \{x_1,\ldots,x_k\} \right \}$, and we use the definition of
$\mathbb{P}_{p,r}^G$ as an image measure. By definition,
\begin{align*}
  \mu_{p,r} ^G & \left [ f(\xi(x_1),\ldots,\xi(x_k))\right ] \\
  &= \sum_{\boldsymbol{v} \in \boldsymbol{V}} \mathbb{P}_{p,r} ^G\left [
    f(\xi(x_1),\ldots,\xi(x_k)) | \{
    \boldsymbol{m}_{\boldsymbol{x}}=\boldsymbol{v} \} \right
  ]\mathbb{P}_{p,r} ^G \left [
    \{\boldsymbol{m}_{\boldsymbol{x}}=\boldsymbol{v} \} \right ]\\
  &= \sum_{\boldsymbol{v} \in \boldsymbol{V}} \nu_p^{\mathcal{E}}
  \otimes \nu_r^{\mathcal{V}} \left [ f(\kappa(v_1),\ldots,\kappa(v_k))
    | \{ \boldsymbol{m}_{\boldsymbol{x}}=\boldsymbol{v} \} \right
  ]\nu_p^{\mathcal{E}}  \left [
    \{\boldsymbol{m}_{\boldsymbol{x}}=\boldsymbol{v} \} \right ]\\
  &= \sum_{\boldsymbol{v} \in \boldsymbol{V}} \nu_r^{\mathcal{V}} \left
    [ f(\kappa(v_1),\ldots,\kappa(v_k))\right ] \nu_p^{\mathcal{E}}
  \left [ \{\boldsymbol{m}_{\boldsymbol{x}}=\boldsymbol{v} \} \right].
\end{align*}
Note that $\nu_r^{\mathcal{V}} \left [
  f(\kappa(v_1),\ldots,\kappa(v_k))\right ]$ is a polynomial in $r$,
so to conclude the proof we only need to prove that for any fixed
$\boldsymbol{x}$ and $\boldsymbol{v}$, $\nu_p ^{\mathcal{E}}\left ( \{
  \boldsymbol{m}(\boldsymbol{x})=\boldsymbol{v} \} \right )$
depends continuously on $p$ on the interval $[0,p_c^G)$.

For $n \geq 1$, write $F_n = \left \{ |C_{x_1}| \leq n, \ldots,
  |C_{x_k}| \leq n \right \}$. It is easy to verify that the event
$\left \{ \boldsymbol{m}_{\boldsymbol{x}}=\boldsymbol{v} \right \}
\cap F_n$ depends on the state of finitely many edges. Hence, $\nu_p
^{\mathcal{E}}\left [ \left \{
    \boldsymbol{m}_{\boldsymbol{x}}=\boldsymbol{v} \right \} \cap F_n
\right ]$ is a polynomial function of $p$.

Fix $p_0 < p_c^G$. For all $p \leq p_0$,
\begin{eqnarray*}
  0 \leq \nu_p^{\mathcal{E}} \left [ \left \{
      \boldsymbol{m}(\boldsymbol{x})=\boldsymbol{v} \right \} \right ] -
  \nu_p ^{\mathcal{E}}\left [ \left \{ \boldsymbol{m}_{\boldsymbol{x}}=\boldsymbol{v}
    \right \} \cap F_n \right ]
  &\leq&   \nu_p^{\mathcal{E}} \left [ F_n^c \right ]\\
  &\leq&   \nu_{p_0}^{\mathcal{E}} \left [ F_n^c \right ]
\end{eqnarray*}
where $\lim\limits_{n \to \infty}\nu_{p_0}^{\mathcal{E}} \left [ F_n^c
\right ]= 0$, since $p_0 < p_c^G$.  So, $\nu_p^{\mathcal{E}} \left [
  \boldsymbol{m}(\boldsymbol{x})=\boldsymbol{v} \right ]$ is a uniform
limit of polynomials on any interval $[0,p_0], \: p_0<p_c^G$, which
implies the desired continuity. \qed

\begin{Remark}
  In the proof we can see that, for fixed $p<p_c^G$, $\mu_{p,r}^G(A)$ is a
  polynomial in $r$.
\end{Remark}

\begin{Remark}
  If $G$ is a graph with uniqueness of the infinite bond cluster in the supercritical regime, then it is possible to verify that $\nu_p^{\mathcal{E}} \left [ \left \{ \boldsymbol{m}(\boldsymbol{x})=\boldsymbol{v} \right \} \right ] $ is continuous in $p$ on the whole interval $[0,1]$. In this case, the continuity given by the Proposition~\ref{cylcont} can be extended to the whole square $[0,1]^2$.
\end{Remark}

\paragraph{Proof of Theorem~\ref{continp}.}
In order to simplify our notations, we write $\mathbb{P}_{p,r}, \nu
_{p}$, $r_c(p)$, for $\mathbb{P}_{p,r}^{\mathbb{Z}^2}, \nu
_{p}^{\mathcal{E}^2}$ and \smash{$r_c^{\mathbb{Z}^2}(p)$}
respectively.  Fix $p_0\in (0,1/2)$ and $\varepsilon >0$
arbitrarily. We will show that there exists $\delta =\delta
(p_0,\varepsilon )>0$ such that for all $p\in (p_0-\delta ,p_0+\delta
)$,
\begin{equation}\label{minus}
  r_c(p) \geq r_c(p_0)-\varepsilon,
\end{equation}
and
\begin{equation}\label{plus}
  r_c(p) \leq r_c(p_0)+\varepsilon.
\end{equation}
Note that by equation (\ref{rcplusrcstaris1}), for all small enough
choices of $\delta >0$ (such that $0\leq p_0\pm \delta <1/2$),
(\ref{minus}) is equivalent to
\begin{equation}\label{plusstar}
  r_c^*(p) \leq r_c^*(p_0)+\varepsilon.
\end{equation}
Below we will show how to find $\delta _1>0$ such that we have
(\ref{plus}) for all $p\in (p_0-\delta _1,p_0+\delta _1)$. One may then
completely analogously find $\delta _2>0$ such that (\ref{plusstar})
holds for all $p\in (p_0-\delta _2,p_0+\delta _2)$, and take $\delta
=\min (\delta _1,\delta _2)$.

\bigskip

Fix $r=r_c(p_0)+\varepsilon $, and define the event $V_n=\{(\xi ,\eta
)\in \Omega _{\mathbb{Z}^2}\times \Omega _{\mathcal{E}_2}:$ there exists
a vertical crossing of $[0,n]\times [0,3n]$ that is black in $\xi \}$.
By ``vertical crossing,'' we mean a self-avoiding path of vertices in
$[0,n]\times [0,3n]$ with one endpoint in $[0,n]\times \{0\}$, and one
in $[0,n]\times \{3n\}$.  Recall also the definition of $M_n$ in
Theorem~\ref{lexpdecay}. By Lemma 2.10 in \cite{BCM}, there exists a
constant $\gamma >0$ such that the following implication holds for any
$p,a\in [0,1]$ and $L\in \mathbb{N}$:
\begin{equation*}\label{eq:finitecriterion}
  \left .
    \begin{array}[l]{rcl}
      (3L+1)(L+1)\nu _{a}(M_{\lfloor L/3\rfloor }) & \leq & \gamma , \\
  \textrm{and } \mathbb{P}_{p,a}(V_L) & \geq & 1-\gamma
    \end{array}
  \right \} \Rightarrow a\geq r_c(p).
\end{equation*}
As usual, $\lfloor x
\rfloor $ for $x>0$ denotes the largest integer $m$ such that $m\leq
x$.  Fix such a $\gamma $.

By Theorem~\ref{lexpdecay}, there exists $N\in \mathbb{N}$ such that
\begin{equation*}\label{p0smallbc}
  (3n+1)(n+1)\nu _{p_0}(M_{\lfloor n/3\rfloor })<\gamma
\end{equation*}
for all $n\geq N$.  On the other hand, since $r>r_c(p_0)$, it follows
from Lemma 2.11 in \cite{BCM} that there exists $L\geq N$ such that
\begin{equation*}\label{p0bigcr}
  \mathbb{P}_{p_0,r}(V_L) > 1-\gamma.
\end{equation*}
Note that both $(3L+1)(L+1)\nu _{p}(M_{\lfloor L/3\rfloor })$ and
$\mathbb{P}_{p,r}(V_L)$ are continuous in $p$ at $p_0$.  Indeed, the
former is simply a polynomial in $p$, while the continuity of the latter
follows from Proposition~\ref{cylcont}.  Therefore, there exists $\delta
_1>0$ such that for all $p\in (p_0-\delta _1,p_0+\delta _1)$,
\begin{eqnarray*}
  (3L+1)(L+1)\nu _{p}(M_{\lfloor L/3\rfloor }) & \leq & \gamma , \\
  \textrm{and } \mathbb{P}_{p,r}(V_L) & \geq & 1-\gamma.
\end{eqnarray*}
By the choice of $\gamma $, this implies that $r\geq r_c(p)$ for all
such $p$, which is precisely what we wanted to prove.

Finding $\delta _2>0$ such that (\ref{plusstar}) holds for all $p\in
(p_0-\delta _2,p_0+\delta _2)$ is analogous: one only needs to substitute
$r_c(p_0)$ by $r_c^*(p_0)$ and ``crossing'' by ``$*$-crossing,'' and the
exact same argument as above works.  It follows that $\delta =\min
(\delta _1,\delta _2)>0$ is a constant such that both (\ref{plus}) and
(\ref{plusstar}) hold for all $p\in (p_0-\delta ,p_0+\delta )$,
completing the proof of continuity on $(0,1/2)$.  Right-continuity at
$0$ may be proved analogously; alternatively, it follows from
Proposition~\ref{contin0}. \qed

\begin{Remark}
  It follows from Theorem~\ref{continp} and equation
  (\ref{rcplusrcstaris1}) that $r_c^*(p)$ is also continuous in $p$ on
  $[0,1/2)$.
\end{Remark}

\section{The critical value functions of tree-like graphs}\label{tree-like-section}

In this section, we will study the critical value functions of
graphs that are constructed by replacing edges of an infinite tree
by a sequence of finite graphs.
We will then use several such constructions in the proofs of our main
results in Section~\ref{proofs-section}.

Let us fix an arbitrary sequence $D_n=(\mathcal{V}_n,\mathcal{E}_n)$
of finite connected graphs and,
for every $n\in \mathbb{N}$, two distinct vertices $a_n,b_n \in \mathcal{V}_n$.
Let $\mathbb{T}_3=(V_3,E_3)$ denote the (infinite) regular tree
of degree $3$, and fix an arbitrary vertex $\rho \in V_3$. Then,
for each edge $e\in E_3$, we denote the end-vertex of $e$ which is closer
to $\rho $ by $f(e)$, and the other end-vertex by $s(e)$.  Let
$\Gamma_D=(\tilde{V},\tilde{E})$ be the graph obtained by replacing
every edge $e$ of $\Gamma _3$ between levels $n-1$ and $n$
(\emph{i.e.}, such that $dist(s(e),\rho )=n$)
by a copy $D_e$ of $D_n$, with $a_n$ and $b_n$
replacing respectively $f(e)$ and $s(e)$. Each vertex $v \in V_3$ is
replaced by a new vertex in $\tilde{V}$, which we denote by $\tilde{v}$.
It is well known that $p_c^{\Gamma _3}=r_c^{\Gamma _3}(0)=1/2$.  Using
this fact and the tree-like structure of $\Gamma _D$, we will be able to
determine bounds for $p_c^{\Gamma _D}$ and $r_c^{\Gamma _D}(p)$.

First, we define $h^{D_n}(p)=\nu_{p}^{\mathcal{E}_n}(a_n \text{ and } b_n \text{ are in the same
  bond cluster})$, and prove the following, intuitively clear, lemma.

\begin{Lemma}
  \label{LemmaBond}
 For any $p\in [0,1]$, the following implications hold:
 \begin{itemize}
 \item[a)] if $\limsup _{n\to \infty }h^{D_n}(p)<1/2$, then $p\leq p_c^{\Gamma_D}$;
 \item[b)] if $\liminf _{n\to \infty }h^{D_n}(p)>1/2$, then $p\geq p_c^{\Gamma_D}$.
 \end{itemize}
\end{Lemma}

\paragraph{Proof.}

We couple Bernoulli bond percolation with parameter $p$ on $\Gamma_D$
with inhomogeneous Bernoulli bond percolation with parameters
$h^{D_n}(p)$ on $\mathbb{T}_3$, as follows. Let $\eta$ be a random
variable with law $\nu_p^{\tilde{E}}$, and define, for each edge $e\in
E_3$, $W(e)=1$ if $\tilde{f(e)}$ and $\tilde{s(e)}$ are connected by a
path consisting of edges that are open in $\eta $, and $W(e)=0$
otherwise. The tree-like structure of $\Gamma _D$ implies that $W(e)$
depends only on the state of the edges in $D_e$, and it is clear that
if $dist(s(e),\rho )= n$, then $W(e)=1$ with probability $h^{D_n}(p)$.

It is easy to verify that there exists an infinite open self-avoiding
path on $\Gamma_D$ from $\tilde{\rho }$ in the configuration $\eta$ if
and only if there exists an infinite open self-avoiding path on
$\mathbb{T}_3$ from $\rho $ in the configuration $W$.  Now, if we
assume $\limsup _{n\to \infty }h^{D_n}(p)<1/2$, then there exists
$t<1/2$ and $N\in \mathbb{N}$ such that for all $n\geq N$,
$h^{D_n}(p)\leq t$.  Therefore, the distribution of the restriction of
$W$ on $L= \{e\in E_3: dist(s(e),\rho )\geq N\}$ is stochastically
dominated by the projection of $\nu _t ^{E_3}$ on $L$. This implies
that, a.s., there exists no infinite self-avoiding path in $W$, whence
$p\leq p_c^{\Gamma_D}$ by the observation at the beginning of this
paragraph. The proof of b) is analogous. \qed

\medskip

We now turn to the DaC model on $\Gamma_D$. Recall that for a vertex
$v$, $C_v$ denotes the vertex set of the bond cluster of $v$. Let
$E_{a_n,b_n}\subset \Omega_{\mathcal{E}_n} \times
\Omega_{\mathcal{V}_n}$ denote the event that $a_n$ and $b_n$ are in
the same bond cluster, or $a_n$ and $b_n$ lie in two different bond
clusters, but there exists a vertex $v$ at distance $1$ from $C_{a_n}$
which is connected to $b_n$ by a black path (which also includes that
$\xi (v)=\xi (b_n)=1$). This is the same as saying that $C_{a_n}$ is
\emph{pivotal} for the event that there is a black path between $a_n$
and $b_n$, \emph{i.e.}, that such a path exists if and only if
$C_{a_n}$ is black. It is important to note that $E_{a_n,b_n}$ is
independent of the color of $a_n$. Define
$f^{D_n}(p,r)=\mathbb{P}_{p,r}^{D_n}(E_{a_n,b_n})$, and note also
that, for $r>0$, $f^{D_n}(p,r)=\mathbb{P}_{p,r}^{D_n}($there is a
black path from $a_n$ to $b_n\mid \xi (a_n)=1)$.

\begin{Lemma}\label{LemmaDaC}
For any $p,r\in [0,1]$, we have the following:
\begin{itemize}
\item[a)] if $\limsup _{n\to \infty }f^{D_n}(p,r)<1/2$, then $r\leq r_c^{\Gamma_D}(p)$;
\item[b)] if $\liminf _{n\to \infty }f^{D_n}(p,r)>1/2$, then $r\geq r_c^{\Gamma_D}(p)$.
\end{itemize}
\end{Lemma}

\paragraph{Proof.}We couple here the DaC model on $\Gamma_D$ with
inhomogeneous Bernoulli \emph{site} percolation on $\mathbb{T}_3$. For
each $v \in V_3\setminus \{ \rho \}$, there is a unique edge $e \in
E_3$ such that $v=s(e)$. Here we denote $D_e$ (\emph{i.e.}, the
subgraph of $\Gamma_D$ replacing the edge $e$) by $D_{\tilde{v}}$, and
the analogous event of $E_{a_n,b_n}$ for the graph $D_{\tilde{v}}$ by
$E_{\tilde{v}}$. Let $(\eta,\xi)$ with values in $\Omega_{\tilde{E}}
\times \Omega_{\tilde{V}}$ be a random variable with law $\mathbb{P}
_{p,r}^{\Gamma _D}$. We define a random variable $X$ with values in
$\Omega_{V_3}$, as follows:
\begin{equation*}
  X(v)=\begin{cases}
    \xi(\tilde{\rho })\!\! & \text{if } v= \rho ,\\
    1 & \text{if the event $E_{\tilde{v}}$ is realized by the
      restriction of $(\eta,\xi)$ to $D_{\tilde{v}}$,} \\
    0 & \text{otherwise.}
  \end{cases}
\end{equation*}

As noted after the proof of Lemma~\ref{LemmaBond}, if $u=f(\left<
  u,v\right>)$, the event $E_{\tilde{v}}$ is independent of the color of
$\tilde{u}$, whence $( E_{\tilde{v}} )_{v \in V_3\setminus \{ \rho \}
}$ are independent. Therefore, as $X(\rho )=1$ with probability
$r$, and $X(v)=1$ is realized with probability $f^{D_n}(p,r)$
for $v\in V_3$ with $dist(v,\rho )=n$ for some $n\in \mathbb{N}$,
$X$ is inhomogeneous Bernoulli site percolation on $\mathbb{T}_3$.

Our reason for defining $X$ is the
following property: it holds for all $v \in V_3\setminus \{ \rho \}$
that
\begin{equation}\label{eq:Connection0}
  \tilde{\rho } \overset{\xi}{\leftrightarrow} \tilde{v} \quad \text{if
    and only if} \quad \rho \overset{X}{\leftrightarrow} v,
\end{equation}
where $x \overset{Z}{\leftrightarrow} y$ denotes that $x$ and $y$ are
in the same \emph{black} cluster in the configuration $Z$. Indeed,
assuming $\tilde{\rho } \overset{\xi}{\leftrightarrow} \tilde{v}$,
there exists a path $\rho =x_0,x_1,\cdots,x_k=v$ in $\Gamma _3$ such
that, for all $0\leq i<k$, $\tilde{x_i} \overset{\xi}{\leftrightarrow}
\tilde{x_{i+1}}$ holds. This implies that $\xi(\tilde{\rho })=1$ and
that all the events $( E_{\tilde{x_i}} )_{0<i\leq k}$ occur, whence
$X(x_i)=1$ for $i=0, \ldots ,k$, so $\rho \overset{X}{\leftrightarrow}
v$ is realized. The proof of the other implication is similar. It
follows in particular from (\ref{eq:Connection0}) that $\tilde{\rho }$
lies in an infinite black cluster in the configuration $\xi$ if and
only if $\rho $ lies in an infinite black cluster in the configuration
$X$.

Lemma~\ref{LemmaDaC} presents two scenarios when it is easy to
determine (via a stochastic comparison) whether the latter event has
positive probability. For example, if we assume that
\mbox{$\liminf_{n\to \infty }f^{D_n}(p,r)>1/2$}, then there exists
$t>1/2$ and $N\in \mathbb{N}$ such that for all $n\geq N$,
$f^{D_n}(p,r)\geq t$. In this case, the distribution of the
restriction of $X$ on $K= \{v\in V_3: dist(v,\rho )\geq N\}$ is
stochastically larger than the projection of $\nu _t ^{E_3}$ on $K$.
Let us further assume that $r>0$. In that case, $X(\rho )=1$ with
positive probability, and $f^{D_n}(p,r)>0$ for every $n\in
\mathbb{N}$. Therefore, under the assumptions $\liminf _{n\to \infty
}f^{D_n}(p,r)>1/2$ and $r>0$, $\rho $ is in an infinite black cluster
in $X$ (and, hence, $\tilde{\rho }$ is in an infinite black cluster in
$\xi$) with positive probability, which can only happen if $r\geq
r_c^{\Gamma_D}(p)$. On the other hand, if $\liminf _{n\to \infty
}f^{D_n}(p,0)>1/2$, then it is clear that $\liminf _{n\to \infty
}f^{D_n}(p,r)>1/2$ (whence $r\geq r_c^{\Gamma_D}(p)$) for all $r>0$,
which implies that $r_c^{\Gamma_D}(p)=0$. The proof of part a) is
similar. \qed

\section{Counterexamples}\label{proofs-section}

In this section, we study two particular graph families and obtain
examples of non-monotonicity and non-continuity of the critical value
function.

\subsection{Non-monotonicity}\label{non-monotonicity-section}

The results in Section~\ref{tree-like-section} enable us
to prove that (a small modification of) the construction
considered by H\"aggstr\"om in the proof of Theorem 2.9 in \cite{HaggstromDaC}
is a graph whose critical coloring value is non-monotone in
the subcritical phase.

\paragraph{Proof of Proposition~\ref{PropNonMononicity}.}

\begin{figure}[ht]
  \begin{center}
    \includegraphics{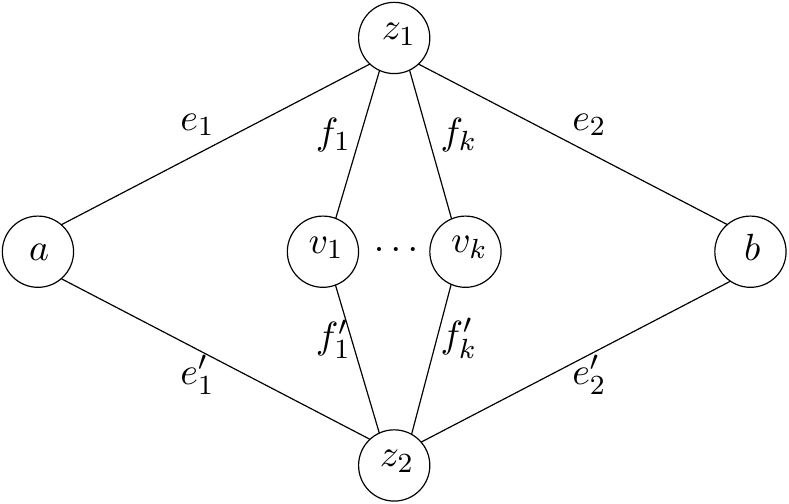}
    \caption{The graph $D^k$.}
    \label{graphDk}
  \end{center}
\end{figure}

Define for $k\in \mathbb{N}$, $D^k$ to be the complete bipartite graph
with the vertex set partitioned into $\{z_1,z_2\}$ and
$\{a,b,v_1,v_2,\ldots ,v_k\}$ (see Figure~\ref{graphDk}). We call
$e_1, e_1'$ and $e_2, e_2'$ the edges incident to $a$ and $b$
respectively, and for $i=1,\ldots,k$, $f_i,f_i'$ the edges incident to
$v_i$. Consider $\Gamma_k$ the quasi-transitive graph obtained by
replacing each edge of the tree $\mathbb{T}_3$ by a copy of $D_k$.
$\Gamma_k$ can be seen as the tree-like graph resulting from the
construction described at beginning of the section, when we start with
the constant sequence $(D_n,a_n,b_n)=(D^k,a,b)$.

We will show below that
it holds for all $k\in \mathbb{N}$ that
\begin{align}
    &p_c^{\Gamma_k}>1/3, \label{fact1}\\
    &r_c^{\Gamma_k}(0)<2/3, \quad \text{and} \label{fact2}\\
    &r_c^{\Gamma_k}(1/3)<2/3. \label{fact3}
\end{align}
Furthermore, there exists $k\in \mathbb{N}$ and $p_0 \in (0,1/3)$ such that
\begin{equation}\label{fact4}
  r_c^{\Gamma_k}(p_0)>2/3.
\end{equation}
Proving (\ref{fact1})--(\ref{fact4}) will finish the proof of
Proposition~\ref{PropNonMononicity} since these inequalities imply
that the quasi-transitive graph $\Gamma_k$ has a non-monotone critical
value function in the subcritical regime.

Throughout this proof, we will omit superscripts in the notation
when no confusion is possible.
For the proof of (\ref{fact1}), recall that $h^{D^k}$ is strictly
increasing in $p$, and $h^{D^k}(p_{D^k})=1/2$.  Since $1-h^{D^k}(p)$ is
the $\nu_{p}$-probability of $a$ and $b$ being in two different bond
clusters, we have that
\begin{equation*}
  1-h^{D^k}(1/3) \geq \nu_{1/3}(\{\text{$e_1$ and $e_1'$ are closed} \}
  \cup \{ \text{$e_2$ and $e_2'$ are closed} \} ).
\end{equation*}
From this, we get that $h^{D^k}(1/3) \leq 25/81$, which proves
(\ref{fact1}).

To get (\ref{fact2}), we need to remember that for fixed $p<p_{D^k}$,
$f^{D^k}(p,r)$ is strictly increasing in $r$, and
$f^{D^k}(p,r_{D^k}(p))=1/2$. One then easily computes that
$f(0,2/3)=16/27>1/2$, whence (\ref{fact2}) follows from
Lemma~\ref{LemmaDaC}.

Now, define $A$ to be the event that at least one edge out of $e_1$,
$e_1'$, $e_2$ and $e_2'$ is open. Then
\begin{eqnarray*}
  f^{D^k}(1/3,2/3) & \geq & \mathbb{P}_{1/3,2/3}(E_{a,b} \mid A)
  \mathbb{P}_{1/3,2/3}(A)\\
  & \geq & \mathbb{P}_{1/3,2/3}(C_b \textrm{ black} \mid A)\cdot 65/81,
\end{eqnarray*}
which gives that $f^{D^k}(1/3,2/3) \geq 130/243 >1/2$, and implies
(\ref{fact3}) by~\ref{LemmaDaC}.

To prove (\ref{fact4}), we consider $B_k$ to be the event that $e_1$,
$e_1'$, $e_2$ and $e_2'$ are all closed and that there exists $i$ such
that $f_i$ and $f_i'$ are both open. One can easily compute that
\[\mathbb{P}_{p,r}(B_k)={(1-p)}^4 \left(1-{(1-p^2)}^k\right),\]
which implies that we can choose $p_0 \in (0,1/3)$ (small) and $k\in
\mathbb{N}$ (large) such that $\mathbb{P}_{p_0,r}(B_{k})>17/18$.
Then,
\begin{eqnarray*}
  f^{D^k}(p_0,2/3) & = &  \mathbb{P}_{p_0,r}(E_{a,b}\mid
  B_k)\mathbb{P}_{p_0,r}(B_k)+\mathbb{P}_{p_0,r}(E_{a,b}\mid
  B_k^c)(1\!-\!\mathbb{P}_{p_0,r}(B_k))\\
  & < & (2/3)^2\cdot 1 + 1\cdot 1/18 ( = 1/2),
\end{eqnarray*}
whence inequality (\ref{fact4}) follows with these choices from
Lemma~\ref{LemmaDaC}, completing the proof. \qed

\subsection{Graphs with discontinuous critical value
  functions}\label{discontgraphs-section}

\paragraph{Proof of Proposition~\ref{nonbounded}.}

\begin{figure}[ht]
  \begin{center}
    \includegraphics{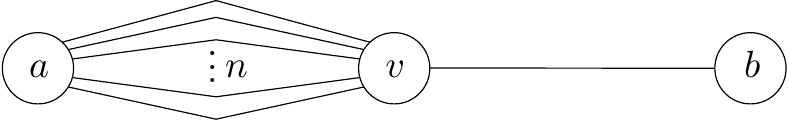}
    \caption{The graph $D_n$.}
    \label{Gn}
  \end{center}
\end{figure}

For $n\in \mathbb{N}$, let $D_n$ be the graph depicted in Figure~\ref{Gn},
and let $G$ be $\Gamma _D$ constructed with this sequence of graphs
as described at the beginning of Section~\ref{tree-like-section}.

It is elementary that $\lim _{n\to \infty }h^{D_n}(p)=p$, whence $p_c^G=1/2$
follows from Lemma~\ref{LemmaBond}, thus $p=0$ is subcritical.
Since $\lim _{n\to \infty }f^{D_n}(0,r)=r^2$, Lemma~\ref{LemmaDaC} gives that
$r_c^G(0)=1/\sqrt{2}$. On the other hand, $\lim _{n\to \infty }f^{D_n}(p,r)=p+(1-p)r$
for all $p>0$, which implies by Lemma~\ref{LemmaDaC} that for $p\leq 1/2$,
\[ r_c^G(p)=\frac{1/2-p}{1-p}\to 1/2 \]
as $p\to 0$, so $r_c^G$ is indeed discontinuous at $0<p_c^G$. \qed

\bigskip

In the rest of this section, for vertices $v$ and $w$, we will write $v\leftrightarrow w$ to denote that
there exists a path of open edges between $v$ and $w$.
Our proof  of Theorem~\ref{boundeddegreeconstruction} will be based on the
Lemma 2.1 in \cite{PSS}, that we rewrite here:

\begin{Lemma}\label{sharpthlemma}
There exists a sequence $G_n=(V^n,E^n)$ of graphs and
$x_n,y_n\in V^n$ of vertices ($n\in \mathbb{N}$) such that
\begin{enumerate}
\item $\nu _{1/2}^{E^n}(x_n\leftrightarrow y_n)>\frac{2}{3}$ for all $n$;
\item $\lim _{n\to \infty }\nu _{p}^{E^n}(x_n\leftrightarrow y_n)=0$ for all $p<1/2$, and
\item there exists $\Delta<\infty $ such that, for all $n$, $G_n$ has
  degree at most $\Delta $.
\end{enumerate}
\end{Lemma}

\bigskip

Lemma~\ref{sharpthlemma} provides a sequence of bounded degree graphs that exhibit
sharp threshold-type behavior at $1/2$. We will use such a sequence as a building block
to obtain discontinuity at $1/2$ in the critical value function in the DaC model.

\paragraph{Proof of Theorem~\ref{boundeddegreeconstruction}.}
We first prove the theorem in the case $p_0=1/2$. Consider the graph
$G_n=(V^n,E^n), x_n,y_n$ $(n\in \mathbb{N})$ as in
Lemma~\ref{sharpthlemma}. We construct $D_n$ from $G_n$ by adding to
it one extra vertex $a_n$ and one edge $\{a_n,x_n\}$. More precisely
$D_n$ has vertex set $V^n\cup\{a_n\}$ and edge set
$E^n\cup\{a_n,x_n\}$.  Set $b_n=y_n$ and let $G$ be the graph $\Gamma
_D$ defined with the sequence $(D_n,a_n,b_n)$ as in
Section~\ref{tree-like-section}.

We will show below that there exists $r_0>r_1$ such that the graph $G$
verify the following three properties:
\begin{enumerate}
\item[(i)] $1/2<p_c^G$
\item[(ii)] $r_c^G(p)\geq r_0$ for all $p<1/2$.
\item[(iii)] $r_c^G(1/2)\leq r_1$.
\end{enumerate}
It implies a discontinuity of $r_c^G$ at $1/2<p_c^G$, finishing the
proof.

One can easily compute $h^{D_n}(p)=p\nu_{p}^{E^n}(x_n\leftrightarrow y_n)$.
Since the graph $G_n$ has degree at most $\Delta$ and the two vertices
$x_n, y_n$ are disjoint, the probability
$\nu_{p}^{E^n}(x_n\leftrightarrow y_n)$ cannot exceed $1-(1-p)^\Delta$. This bound
guarantees the existence of $p_0>1/2$ independent of $n$ such that $h^{D_n}(p_0)<1/2$ for all $n$, whence
Lemma~\ref{LemmaBond} implies that $1/2<p_0\leq p_c^G$.

For all $p\in [0,1]$, we have
\begin{equation*}\label{fforG}
  f^{D_n}(p,r)\leq \left(p+r(1-p)\right)\left(\nu _{p}^{E^n}(x_n\leftrightarrow y_n)+r (1-\nu _{p}^{E^n}(x_n\leftrightarrow y_n))\right).
\end{equation*}
which gives that
$\underset{n\to\infty}{\lim}f^{D_n}(p,r)<\left(\frac{r+1}{2}\right)\cdot
r$. Writing $r_0$ the positive solution of $r(1+r)=1$, we get that
$\underset{n\to\infty}{\lim}f^{D_n}(p,r_0)<1/2$ for all $p<1/2$, which implies by Lemma~\ref{LemmaDaC} that
$r_c^G(p)\geq r_0$. 

On the other hand, $f^{D_n}(1/2,r)\geq \nu
_{p}^{E^n}(x_n\leftrightarrow y_n)\left(\frac{1+r}{2}\right)$, which gives by
Lemma~\ref{sharpthlemma} that $\underset{n\to\infty}{\lim}
f^{D_n}(1/2,r)>\frac{2}{3}\cdot\frac{1+r}{2}$. Writing $r_1$ such that
$\frac{2}{3}(1+r_1)=1$, it is elementary to check that $r_1<r_0$ and
that $\underset{n\to\infty}{\lim}
f^{D_n}(1/2,r_1)>1/2$.
Then, using Lemma~\ref{LemmaDaC}, we conclude that
$r_c(1/2)\leq r_1$. \qed

\paragraph{Acknowledgments.} We thank Jeff Steif for suggesting (a
variant of) the graph that appears in the proof of
Theorem~\ref{boundeddegreeconstruction}. V.B. and V.T. were supported
by ANR grant 2010-BLAN-0123-01.

\end{document}